\documentclass[a4paper] {article}
\usepackage[latin1]{inputenc}
\usepackage[english]{babel}
\usepackage{graphicx}
\usepackage{amsmath,amsfonts,amssymb,amsthm,amsopn,amscd}
\usepackage{bbm,wasysym,stmaryrd}
\usepackage{subfigure}
\usepackage{apalike}

\newcommand{\der}[2]{\frac{\text{d} #1}{\text{d} #2}}	        

\newtheorem{Assump}{Assumption}

\title{Sensitivity to the cutoff value in the quadratic adaptive integrate-and-fire model}
\author{Jonathan Touboul}

\begin{document}
 \maketitle
 \section*{Abstract}
 {\bf The quadratic adaptive integrate-and-fire model \cite{izhikevich:03,izhikevich:07} is recognized as very interesting for its computational efficiency and its ability to reproduce many behaviors observed in cortical neurons. For this reason it is currently widely used, in particular for large scale simulations of neural networks. This model emulates the dynamics of the membrane potential of a neuron together with an adaptation variable. The subthreshold dynamics is governed by a two-parameter differential equation, and a spike is emitted when the membrane potential variable reaches a given cutoff value. Subsequently the membrane potential is reset, and the adaptation variable is added a fixed value called the spike-triggered adaptation parameter. We show in this note that when the system does not converge to an equilibrium point, both variables of the subthreshold dynamical system blow up in finite time whatever the parameters of the dynamics. The cutoff is therefore essential for the model to be well defined and simulated. The divergence of the adaptation variable makes the system very sensitive to the cutoff: changing this parameter dramatically changes the spike patterns produced. Furthermore from a computational viewpoint, the fact that the adaptation variable blows up and the very sharp slope it has when the spike is emitted implies that the time step of the numerical simulation needs to be very small (or adaptive) in order to catch an accurate value of the adaptation at the time of the spike. It is not the case for the similar quartic \cite{touboul:08b} and exponential \cite{brette-gerstner:05} models whose adaptation variable does not blow up in finite time, and which are therefore very robust to changes in the cutoff value. }

\section{Introduction}
During the past few years, in the neuro-computing community, the problem of finding a computationally simple and biologically realistic model of neuron has been widely studied, in order to be able to compare experimental recordings with numerical simulations of large-scale brain models. The key problem is to find a model of neuron realizing a compromise between its simulation efficiency and its ability to reproduce what is observed at the cell level, often considering in-vitro experiments \cite{koch-segev:98,izhikevich:04,rinzel-ermentrout:89}. Among the variety of computational neuron models, nonlinear spiking models with adaptation have recently been studied by several authors
\cite{izhikevich:04,brette-gerstner:05,touboul:08b} and seem to stand out. They are relatively simple, i.e. mathematically tractable, efficiently implemented, and able to reproduce a large number of electrophysiological signatures such as bursting or regular spiking. These models satisfy the equations:
\begin{equation}\label{eq:GeneralModel}
 \begin{cases}
  \der{v}{t} &= F(v)-w+I \\
  \der{w}{t} & = a\,(b\,v -w)
 \end{cases}
\end{equation}
where $a$ and $b$ are non-negative parameters and $F(v)$ is a regular strictly convex function satisfying assumption:
\begin{Assump}\label{Assump:Blow}
 There exists $\varepsilon>0$ and $\alpha>0$ for which $F(v)\geq \alpha v^{1+\varepsilon}$ when $v\to \infty$ (we will say that $F$ grows faster than $v^{1+\varepsilon}$ when $v \to \infty$).
\end{Assump}
A spike is emitted at the time $t^*$ when the membrane potential $v$ reaches a cutoff value $\theta$. At this time, the membrane potential is reset to a constant value $c$ and the adaptation variable is updated to $w(t^*)+ d$ where $w(t^*)$ is the value of the adaptation variable at the time of the spike and $d>0$ is the spike-triggered adaptation parameter. 

For these models we prove in section \ref{sec:divergence} that the membrane potential blows up in finite time. Among these models, the \emph{quadratic adaptive} model \cite{izhikevich:04} corresponds to the case where $F(v)=v^2$, and has been recently used by Eugene Izhikevich and coworkers \cite{izhikevich-edelman:08} in very large scale simulations of neural networks. The \emph{adaptive exponential} model \cite{brette-gerstner:05} corresponds to the case where $F(v)=e^v$, has the interest that its parameters can be related to electrophysiological quantities, and has been successfully fit to intracellular recordings of pyramidal cells \cite{clopath-jolivet-etal:07,jolivet-kobayashi-etal:08}. The \emph{quartic} model \cite{touboul:08b} corresponds to the case where $F(v)=v^4+2a\,v$ and has the advantage to of being able to reproduce all the behaviors featured by the other two and also self-sustained subthreshold oscillations which are of particular interest to model certain nerve cells.

In these models, the reset mechanism makes critical the value of the adaptation variable at the time of the spike. Indeed, when a spike is emitted at time $t^*$, the new initial condition of the system \eqref{eq:GeneralModel} is $(c,w(t^*)+d)$. Therefore, this value governs the subsequent evolution of the membrane potential, and hence the spike pattern produced. For instance in \cite{touboul-brette:08,touboul-brette:08c}, the authors show that the sequence of reset locations after each spike time shapes the spiking signature of the neuron. 

Hence characterizing the reset location of the adaptation variable is essential to characterize the spiking properties of these models. To this end, we precisely study in this note the orbits of equation \eqref{eq:GeneralModel} in the phase plane $(v,w)$ in order to characterize the value of the adaptation variable at the time of the spike. We prove in section \ref{sec:divergence} that the adaptation variable diverges when $v\to \infty$ in the case of the quadratic model and converges in the cases of the exponential and of the quartic model, and study in section \ref{sec:discussion} the consequences of this fact on the spiking signatures and on numerical simulation methods.

\section{Adaptation variable at the times of the spikes}\label{sec:divergence}
As we can see in equation \eqref{eq:GeneralModel}, the greater the membrane potential the greater the derivative of the adaptation variable. When the membrane potential blows up, the adaptation variable may either remain bounded or blow up, depending on the shape of the divergence of $v$. When this divergence is not fast enough, the adaptation variable simultaneously blows up. 

We prove here that for the models satisfying assumption \ref{Assump:Blow} the membrane potential blows up in finite time. We also prove that for quadratic adaptive model\footnote{We can prove more generally that when $F(v)/v^2$ tends to a finite constant (possibly $0$), the adaptation variable will blow up when the membrane potential blows up} the adaptation variable blows up at the same time as a logarithmic function of $v$, whereas if there exists $\varepsilon>0$ such that $F(v)$ grows faster than $v^{2+\varepsilon}$ when $v \to \infty$, then the adaptation variable remains bounded when $v \to \infty$. 
 
In \cite{touboul:08b}, we have seen that there exists possibly one stable fixed point for system \ref{eq:GeneralModel}, which corresponds to a resting state. In \cite{touboul-brette:08c}, we prove that all the orbits of the system that do not converge to this stable fixed point will be trapped after a finite time in a zone fully included in the half space $\{w<b\,v\}$ called the \emph{spiking zone}\footnote{In the case where the subthreshold system has no fixed point this property can be derived from the shape of the vector field in the phase plane, as well as in the case where the initial condition $(v,w)$ is such that $v$ is greater than the largest $v$-value of the fixed points (the biggest solution of $F(v)-b\,v+I=0$) and $w\leq b\,v$: in this case the vector field on the line $w=b\,v$ implies that the trajectory keeps trapped in this zone. In the case where there exist fixed points, the proof is slightly more complex and involves the description of the stable manifold of the saddle fixed point.}. Denote $t_0$ a time such that the orbit is inside the spiking zone. In this zone, we have 
\[\der{v}{t} \geq F(v) - b\,v + I \]
It is simple to prove that the solution of the equation
\[\begin{cases}
   \der{u}{t} = F(u) - b\,u + I  \\
   u(t_0)=v(t_0)
  \end{cases}\]
blows up in finite time under the assumption \ref{Assump:Blow}\footnote{For the quadratic model we can get analytic expressions of the solutions involving the tangent function, and therefore can derive an upperbound of the explosion time.}. Using Gronwall's theorem \cite{gronwall:19} we conclude that $v(t)\geq u(t)$ and hence $v$ blows up in finite time.

To prove the divergence of the adaptation variable when the membrane potential blows up in the case of the quartic model, we study the orbit of a solution $(v(t),w(t))$ of the differential system \eqref{eq:GeneralModel} such that the membrane potential blows up at time $t^*$, and characterize the behavior of $w(t)$ in function of $v(t)$. In the spiking zone, we have seen that $w(t) \leq b\,v(t)$ and therefore $F(v)-w+I \geq F(v)-b\,v+I$ which tends to infinity when $v$ tends to infinity. Since $v(t)$ blows up there exists a time $t_1 \in [t_0,t^*)$ such that we will have $F(v(t))-w(t)+I \geq k >0$ for all $t \in [t_1, t^*)$. We denote $(v_1:= v(t_1),w_1 := w(t_1))$. After time $t_1$, because of this inequality, the trajectory in the phase plane can be written as the graph of a function $W(v)$ that satisfies the equation:
\begin{equation}\label{eq:Trajectoire}
 \begin{cases}
  \der{W}{v} = \frac{a\,(b\,v-W)}{F(v)-W+I}\\
   W(v_1)  = w_1
 \end{cases}
\end{equation}
(i.e. $w(t) = W(v(t))$ for $t\in [t_1,t^*)$). Since $w(t)$ is increasing for $t\in [t_1,t^*)$, we necessarily have:
\begin{equation}\label{eq:InegGronwall}
 \der{W}{v} \geq \frac{a\,(b\,v-W)}{F(v)-w_{1}+I} 
\end{equation}
Therefore Gronwall's theorem \cite{gronwall:19} ensures us that the solution of equation \eqref{eq:Trajectoire} will be lowerbounded for $v\geq v_1$ by the solution of the linear ordinary differential equation:
\begin{equation}\label{eq:Gronwall}
 \begin{cases}
  \der{z}{v} = \frac{a\,(b\,v-z)}{F(v)-w_{1}+I}\\
  z(v_1) = w_1
  \end{cases}
\end{equation}
that reads:
\[z(v) = \left(\int_{v_1}^v \frac{a\,b\,u}{F(u)-w_{1}+I}e^{-g(u)}\, du + w_1 \right) e^{g(v)}\]
where $g(v) = -\int_{v_1}^{v}\frac{a\,du}{F(u)-w_{1}+I}$. Because of assumption \ref{Assump:Blow}, the integrand is integrable, and the function $g$ has a finite limit $g(\infty)$ when $v\to\infty$. The exponential terms will hence converge when $v\to \infty$. But the integral involved in the particular solution diverges in the quadratic case\footnote{or when $F(v)$ grows slower than $v^2$,}, since the integrand is equivalent when $u \to \infty$ to 
\[\frac{a\,b}{u}e^{-g(\infty)}\]

Hence the solution of the linear differential equation \eqref{eq:Gronwall} tends to infinity when $v\to \infty$ faster than a logarithmic function of $v$, and so does $W(v)$, and hence $w(t)$ blows up at the time when $v(t)$ blows up. 

Let us now upperbound the adaptation variable on the orbits of the system. Using the same notations, since $w_1\leq w(t) \leq b\,v(t)$ for $t \in [t_1,t^*)$, we have:
\begin{equation}\label{eq:Majoration}
 \der{W}{v}\leq \frac{a\,(b\,v-w_1)}{F(v)-b\,v + I}
\end{equation}
and hence 
\[W(v)\leq w_1+\int_{v_1}^v \frac{a\,(b\,u-w_1)}{F(u)-b\,u + I}\,du\]

In the case where $F(u)=u^2$ this integral is bounded by a logarithmic function of $v$ and in the case where $F(u)$ grows faster than $u^{2+\varepsilon}$, this integral converges when $v\to \infty$. Furthermore, since $W$ is an increasing upperbounded function, it converges when $v\to \infty$.

We therefore conclude that in the case of the quadratic adaptive model, the adaptation variable blows up at the explosion time of the membrane potential variable $v$ and this divergence is logarithmic in $v$ (see figure \ref{fig:Bifurcs} (h)), and in the case of the quartic and exponential models, the adaptation variable converges. Note that the adaptation variable diverges whatever the parameters of the system provided the membrane potential variable blows up. The value of the adaptation variable at the cutoff $\theta$ is simply given by $W(\theta)$, that depends on the parameters of the system and of the initial condition. In the case of the quadratic model it is an unbounded increasing function of $\theta$ , and in the quartic and exponential models, a converging function of $\theta$. 

\section{Consequences}\label{sec:discussion}
The divergence of the adaptation variable at the times of the spikes significantly impacts the theoretical, qualitative and computational analysis of the model. It appears to be a critical parameter of the quadratic model.
 

We have seen that changing the cutoff value resulted in changing the value of the adaptation variable at the times of the spikes. Let $(v_0,w_0)$ be an initial condition for the system \eqref{eq:GeneralModel}. If the neuron fires, its membrane potential will reach the cutoff value $\theta$ at a given time. Since the membrane potential blows up in finite time, the time of the first spike emitted will not be very sensitive to changes in the cutoff value provided it is high enough. But the after-spike reset location $(c,W(\theta)+d)$ will significantly change when varying $\theta$. The whole subsequent evolution of the system is therefore affected, as soon as the second spike is emitted. Thus the spike pattern produced depends on the cutoff value. 

In the case of the quartic and exponential models, the adaptation variable converges when the cutoff tends to infinity. Therefore, the model defined by \eqref{eq:GeneralModel} with an infinite cutoff value is mathematically well defined. In that case, a spike is emitted when the membrane potential blows up and subsequently we reset the membrane potential to a fixed value $c$ and add to the value of adaptation variable at the explosion time the spike-triggered adaptation parameter. We call this system the \emph{intrinsic} system. The behavior of the system and the spike patterns it produces can be mathematically studied (see \cite{touboul-brette:08,touboul-brette:08c}). Interestingly, these intrinsic spike patterns undergo bifurcations with respect to the parameters of the model. When considering a finite cutoff, the model (or the numerical simulation) will approximate these intrinsic behaviors provided that the cutoff threshold is high enough. The sensitivity to the cutoff in these cases will hence be very limited except in very small regions of the parameter space around the bifurcations of the intrinsic system. Unfortunately, for the quadratic model, no intrinsic behavior can be defined because of the divergence of the adaptation variable: the behaviors it produces will depend on the choice of the threshold.

\begin{figure}
 \begin{center}
 \includegraphics[width=.7\textwidth]{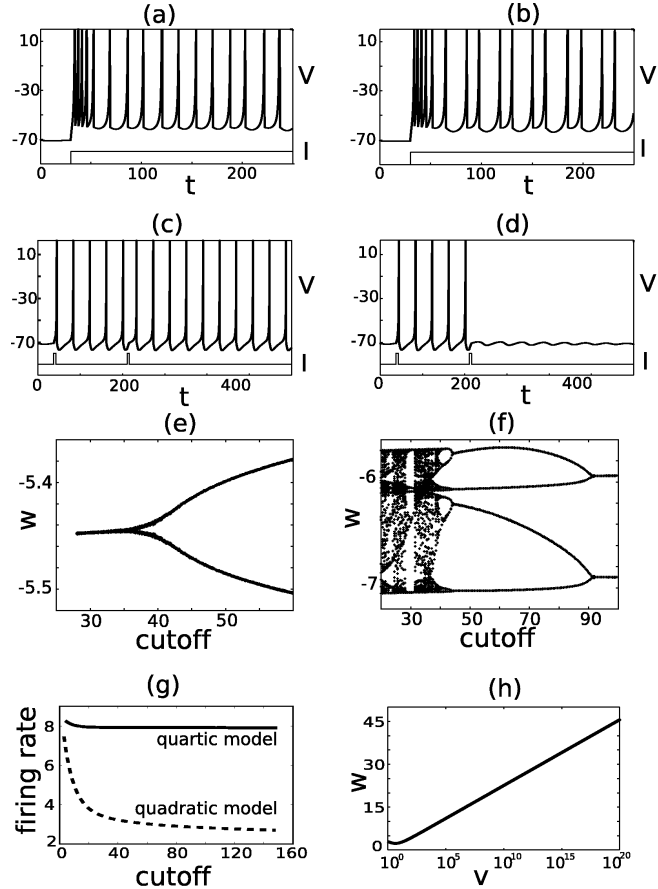}
 \end{center}
 \caption{Sensitivity of the spike patterns with respect to the cutoff value for the quadratic model, for different set of parameters. Parameters used: $(A) = \{a=0.02; b=0.19;  c=-60;  d= 1.419, I=10.25\}$; $(B) = \{a=0.1;  b=0.26; c=-60;  d=0;\}$, $(C) = \{a=0.02, b=0.19, c=-57.7, d=1.15, I=10.377\}$. Simulations done with an Euler scheme with time steps ranging from $10^{-4}$ to $10^{-3}$. For figure (a) and (b) the parameters used are $(A)$ with cutoff of $30$ and $45$ respectively: a small increase of the cutoff results in a sharp transition from spiking to bursting linked with a period doubling bifurcation for the adaptation value at the reset represented in figure (e). Since there is a bifurcation, the qualitative change is very sharp. Figures (c) and (d) corresponds to the parameters $(B)$ with cutoffs value $32.9$ and $33$ respectively. Changing the cutoff results in two very different global behaviors. Fig. (e) and (f) represent the stationary sequence of reset values as functions of the threshold $\theta$. Figure (f) corresponds to the set of parameters $(C)$ for cutoff values ranging from $20$ to $100$: an intricate bifurcation structure appears. Figure (g) shows the convergence of the firing rate to the intrinsic firing rate in the case of the quartic model, while the firing rate of the quadratic model regularly decreases to 0. (h): Divergence of the variable $w$ in function of $v$ in semilogarithm axis.
 }
 \label{fig:Bifurcs}
\end{figure}

First of all, we have seen that the dependency of this reset location in the quadratic model is a logarithmic function of $\theta$, which makes the variations of the reset value in function of the cutoff unbounded but quite slow. Small changes in the cutoff slightly impact the value of the reset adaptation variable. For instance if we consider the firing rate of a neuron in the case where the system has no fixed point, increasing the cutoff value results in the case of the quadratic model in a slow continuous decrease of the firing rate of the neuron that tends to zero as the cutoff increases, whereas the firing rate converges for the quartic model to the related intrinsic firing rate (see figure \ref{fig:Bifurcs}.(g)). 

When considering the spike patterns produced, the effects of changes in the cutoff value for the quadratic model are much more dramatic. Indeed, the sequence of adaptation values at the times of spikes shapes the spike pattern produced: for instance, regular spiking corresponds to the convergence of this sequence, and bursting to cycles in this sequence. These properties are very sensitive to changes in the parameters of the model: bifurcations between different spike patterns, and even chaos appear when the model's parameters vary (see \cite{touboul-brette:08,touboul-brette:08c,naud-macille-etal:08b}). In the case of the quadratic model, we have seen that these adaptation values strongly depend on the cutoff. Therefore, since the dependency on the cutoff is unbounded, from a given initial condition and for fixed values of the parameters, increasing the cutoff may result in crossing many bifurcation lines, and hence in producing many different behaviors. We present in figure \ref{fig:Bifurcs} a graph showing that bifurcations and chaos occur with respect to the cutoff value, in the usual range of simulation parameters. For instance, a period doubling bifurcation appears when varying the cutoff value (in figure \ref{fig:Bifurcs}(e) we give a graph of the stationary reset values in function of the threshold $\theta$), that results in abruptly switching from a regular spiking behavior to a bursting behavior (figures \ref{fig:Bifurcs}(a) and \ref{fig:Bifurcs}(b)). More complex bifurcation structures involving chaotic patterns also appear, and in this case, infinitesimal changes in the cutoff value result in dramatic changes in the behavior. This raises the question of the meaning of the cutoff value in these ranges of parameters (see figure \ref{fig:Bifurcs}(f)). Changing the cutoff in that case makes the system switch between chaotic spiking, bursts with $8$, $4$ and eventually $2$ spikes, for the cutoff values considered. And this behavior will not be observed only for very particular values of the parameters of the system. Depending on the extension of the interval where the cutoff value varies, quite a large set of parameters will present bifurcations in the nature of the emitted spike train.

Because of this sensitivity, the cutoff value and the different parameters of the model have to be very carefully evaluated in order to quantitatively fit datasets. In this context the meaning of the threshold and therefore the problem of its accurate evaluation has to be specifically addressed in the case of the quadratic model, since it has no clear biophysical interpretation.

Eventually, from the numerical viewpoint, the unboundedness of the adaptation variable and of its time derivative at the explosion times of the membrane potential makes the accurate computation of this value very difficult. In particular, the time step necessary to accurately estimate this value has to be very small (or to be adaptive as a function of the value of the membrane potential variable) in order to obtain the right spike pattern. These remarks relativize the statement that this model can be efficiently simulated since very accurate methods have to be implemented in order to correctly evaluate the adaptation variable at the time of the spike.

These remarks do not apply for the models where the adaptation variable converges at the times of the spikes. In these cases, the system has intrinsic properties that make the times of the spike and the adaptation variable at these times robust to the choice of the cutoff value provided it is big enough and the numerical simulations less sensitive to the choice of the time step. However, the exponential integrate-and-fire model involves an exponential function which increases very fast as a function of the membrane potential. This fact implies that the time step chosen for the simulation has to be very small and makes this system less suitable for large scale simulations. In particular, the numerical scheme must be designed carefully in order to capture the value of the adaptation variable at the times of the spike. Simulating the orbit in the phase plane for instance would result in computing very accurately this essential value.

\section*{Conclusion}
In this note we proved that the adaptation variable of the adaptive quadratic model blew up at the times of the spikes whereas it converged for the quartic and the adaptive exponential models. This property has some important implications that are discussed in the note. From a theoretical point of view, we showed that the nature of the spike patterns produced undergoes bifurcations with respect to the cutoff value, and this made the system very sensitive to this parameter: small changes in the value of this parameter can deeply affect the nature of the spiking pattern. From a quantitative viewpoint, it raises the question of how to evaluate this threshold in order to fit datasets, and from a numerical viewpoint, it has implications on the efficiency of the simulation algorithms to use. The convergence of this value for models having a faster blow up at the times of the spike, such as the quartic or the exponential adaptive models, implies that the system presents intrinsic spiking properties which can be mathematically studied, efficiently simulated and robust to changes in the cutoff value.

\section*{Acknowledgments}
The author warmly acknowledges Bruno Cessac, Olivier Faugeras and Fran\c{c}ois Grimbert for interesting discussions and suggestions.

\end{document}